\documentclass{article}

\usepackage{amssymb}
\usepackage{amsfonts}
\usepackage{amsmath}
\usepackage{amsthm}
\usepackage{mathrsfs}
\usepackage{dsfont}
\usepackage{theoremref}

\newtheorem{thm}{Theorem}[section]
\newtheorem{lem}[thm]{Lemma}
\newtheorem{cor}[thm]{Corollary}
\newtheorem{prp}[thm]{Proposition}

\theoremstyle{definition}
\newtheorem{dfn}[thm]{Definition}

\usepackage{a4wide}

\title{
\textbf{\huge{MAD Families of Projections on $l^2$ and Real-Valued Functions on $\omega$}}\\
}
\author{{\Large Tristan Bice}}
\date{}

\begin{document}

\maketitle

\paragraph{Abstract} Two sets are said to be almost disjoint if their intersection is finite.  Almost disjoint subsets of $[\omega]^\omega$ and $\omega^\omega$ have been studied for quite some time.  In particular, the cardinal invariants $\mathfrak{a}$ and $\mathfrak{a}_e$, defined to be the minimum cardinality of a maximal infinite almost disjoint family of $[\omega]^\omega$ and $\omega^\omega$ respectively, are known to be consistently less than $\mathfrak{c}$.  Here we examine analogs for functions in $\mathbb{R}^\omega$ and projections on $l^2$, showing that they too can be consistently less than $\mathfrak{c}$.

\section{Introduction}

\subsection{Projections}

A cardinal invariant is, broadly speaking, any cardinal with some combinatorial definition whose value can not be determined in ZFC.  For example, as mentioned above, there is the cardinal invariant $\mathfrak{a}$ whose value must lie within $\aleph_1$ and $\mathfrak{c}$ but (if CH fails) can consistently take almost any value in between.  Indeed, we can define a general $R_\kappa$-disjointness number as follows.
\begin{dfn}
Let $\kappa$ be a cardinal, $S$ a set and $R$ a (symmetric) binary relation on $S$.  We define \[\mathfrak{a}_\kappa((S,R))=\inf\{|A|:|A|\geq\kappa\wedge A\subseteq S\wedge\forall a,b\in A(aRb)\wedge\forall s\in S\exists a\in A\neg(sRa)\}.\]
\end{dfn}
So $\mathfrak{a}=\mathfrak{a}_{\aleph_0}(([\omega]^\omega,\not|^*))$, where $A \not|^*\ B\Leftrightarrow|A\cap B|<\infty$.  Note here that, as there are trivial examples of finite maximal almost disjoint families of subsets of $\omega$, we need to eliminate these from consideration to define the non-trivial cardinal invariant $\mathfrak{a}$.  However, with the relations $R$ we will be concerned with, there will be no such finite maximal $R$-disjoint families, in which case we may omit the restriction that $|A|\geq\aleph_0$ and accordingly omit the subscript, i.e. $\mathfrak{a}((S,R))=\mathfrak{a}_0((S,R))$.  Also, as is common practice with such definitions, we shall often abbreviate $\mathfrak{a}_\kappa((S,R))$ by $\mathfrak{a}_\kappa(S)$ or $\mathfrak{a}_\kappa(R)$ when the context makes it clear so, for example, $\mathfrak{a}=\mathfrak{a}_{\aleph_0}(\not|^*)$.  

These kinds of cardinal invariants are the focus of this paper.  However, they are certainly not the only kinds of cardinal invariants that can be defined \textendash\, see \cite{a} for a comprehensive introduction to cardinal invariants defined from subsets of $[\omega]^\omega$, for example.  In the past few years, analogous cardinal invariants have been defined from subsets of $\mathcal{P}_\infty(H)$, the collection of infinite rank projections (idempotent self-adjoint operators $P\in\mathcal{B}(H)$, i.e. satisfying $P^2=P$ and $P^*=P$) on a separable infinite dimensional Hilbert space $H$ (which, in fact, is unique up to an isometric isomorphism, $l^2$ being the typical example of such a Hilbert space).  For example in \cite{b} the cardinal invariant $\mathfrak{a}^*$ is defined to be the minimum cardinality of a maximal infinite almost orthogonal subset of $\mathcal{P}_\infty(H)$, according to the following definition.
\begin{dfn}
$P,Q\in\mathcal{P}(H)$ are \emph{almost orthogonal}, written $P\bot^*Q$, if $PQ$ is a compact operator.
\end{dfn}
So, in our notation, $\mathfrak{a}^*=\mathfrak{a}_{\aleph_0}(\bot^*)=\mathfrak{a}_{\aleph_0}((\mathcal{P}_\infty(H),\bot^*))$.

Here we study a slightly different analog of $\mathfrak{a}$ from subsets of $\mathcal{P}_\infty(H)$, namely $\mathfrak{a}(\top^*)$.  Indeed, it is a general fact that cardinal invariants defined using infinite subsets of $\omega$ have (at least) two different natural analogs for infinite rank projections, depending on whether we reinterpret $\not|^*$ as $\bot^*$ or $\top^*$.  To define $\top^*$, recall that the compact operators $\mathcal{K}(H)$ on $H$ form a closed ideal and so we can form the factor C$^*$-algebra, known as the Calkin algebra of $H$, denoted by $\mathcal{C}(H)$, with canonical homomorphism $\pi:\mathcal{B}(H)\mapsto\mathcal{C}(H)$.

\begin{dfn}
$P,Q\in\mathcal{P}(H)$ are \emph{almost disjoint}, written $P\top^*Q$, if $||\pi(PQ)||<1$.
\end{dfn}

Note that $P,Q\in\mathcal{P}(H)$ are almost orthogonal if and only if $||\pi(PQ)||=0$, so $P\bot^*Q$ implies $P\top^*Q$.  In fact, almost orthogonality is equivalent to almost disjointness together with commutivity (modulo compact operators).

\begin{prp}
For $P,Q\in\mathcal{P}_\infty(H)$,
\begin{equation}
P\bot^*Q\quad\Leftrightarrow\quad\pi(PQ)=\pi(QP)\ \wedge\ P\top^*Q.\nonumber
\end{equation}
\end{prp}

\paragraph{Proof:} If $P\bot^*Q$ then $\pi(QP)=\pi((PQ)^*)=\pi(PQ)^*=0^*=0=\pi(PQ)$ and, as mentioned above, $P\top^*Q$.  For the other direction note that $\pi(PQ)=\pi(QP)$ implies that $\pi(PQ)$ is a projection and hence $\pi(PQ)=0$, because $||p||=1$ for all non-zero projections $p\in\mathcal{C}(H)$. $\Box$\\

Indeed, most research on cardinal invariants defined from projections up till now has focused on (modulo compact) commutative subsets of $\mathcal{P}_\infty(H)$.  This is probably because even basic questions about non-commutative projections had not been answered, such as the question of when two non-commutative $P,Q\in\mathcal{P}(H)$ have a greatest $\leq^*$-lower bound, where \[P\leq^*Q\quad\Leftrightarrow\quad\pi(PQ)=\pi(P)\quad\Leftrightarrow\quad PQ-P\textrm{ is compact.}\]

However, recent progress in this direction has been made in \cite{c}, where it is also shown that $||\pi(PQ)||<1$ is equivalent to saying $P$ and $Q$ have no non-trivial (i.e. infinite rank) lower bound with respect to $\leq^*$.  So $\mathfrak{a}(\top^*)$ is the minimum cardinality of a maximal antichain of infinite rank projections with respect to $\leq^*$.  As $\mathfrak{a}$ is the minimum cardinality of a maximal infinite antichain of infinite subsets of $\omega$ with respect to $\subseteq^*$ (inclusion modulo finite subsets), $\mathfrak{a}(\top^*)$ is perhaps the more natural purely order-theoretic analog of $\mathfrak{a}$ and thus more deserving of the notation $\mathfrak{a}^*$.  However, to avoid confusion, we shall avoid this notation entirely and refer only to $\mathfrak{a}(\bot^*)$ and $\mathfrak{a}(\top^*)$.

\subsection{Functions on $\omega$}

There are also structures that lie somewhere in between $([\omega]^\omega,\not|^*)$ and $(\mathcal{P}_\infty(H),\top^*)$.  To see this, consider $H=\bigoplus_\omega\mathbb{F}^2$ (here $\mathbb{F}$ denotes the scalar field of $H$, either the reals $\mathbb{R}$ or the complex numbers $\mathbb{C}$).  Given $(\theta_n)\subseteq[0,\pi/2]$ set $\mathbf{v}_n(m)=\delta_{m,n}(\cos\theta_n,\sin\theta_n)\in H$ and let $P_{(\theta_n)}\in\mathcal{P}(H)$ be such that $\mathcal{R}(P_{(\theta_n)})=V_{(\theta_n)}=\overline{\mathrm{span}}(\mathbf{v}_n)$.  Then, for $(\theta_n),(\phi_n)\subseteq[0,\pi/2]$, \[||\pi(P_{(\theta_n)}P_{(\phi_n)})||=\lim\inf_n\cos(\theta_n-\phi_n),\] so
$P_{(\theta_n)}\top^*P_{(\phi_n)}\Leftrightarrow\lim\inf|\theta_n-\phi_n|>0$.  This inspires the following definition.
\begin{equation}\label{f|g}
f\not\approx g\Leftrightarrow\lim\inf|f(n)-g(n)|>0.
\end{equation}
Of course, $\mathfrak{a}([0,\pi/2]^\omega,\not\approx)$ will only equal the minimum cardinality of a maximal antichain of a certain subcollection of projections, but still provides a good intuitive basis for studying antichains of arbitrary projections.  In fact, more generally, we define the following.
\begin{dfn}\thlabel{f|g2}
$f,g\in\mathscr{P}(\mathbb{R})^\omega$ are \emph{lim-inf disjoint}, written $f\not\approx g$, if there exists $m\in\omega$ and $\epsilon>0$ such that $\forall k>m\forall s\in f(k)\forall t\in g(k)(|s-t|>\epsilon)$.
\end{dfn}

For a function $f$ on $\omega$ define $f'$ on $\omega$ by $f'(n)=\{f(n)\}$, for all $n\in\omega$.  So, for $f,g\in\mathbb{R}^\omega$, we have $f\not\approx g$ according to (\ref{f|g}) if and only $f'\not\approx g'$ according to \thref{f|g2}.  Context should make it clear which definition is being used.  If $\mathcal{F}$ is a collection of functions on $\omega$, let $\mathcal{F}_\infty$ be the subset of all $f\in\mathcal{F}$ such that $\{n:f(n)\neq0\}$ is infinite.  We are primarily interested in $\mathfrak{a}(\mathbb{R}^\omega)$, $\mathfrak{a}(([\mathbb{R}]^{<\omega}\backslash\{0\})^\omega)$, $\mathfrak{a}(([\mathbb{R}]^{\leq1})^\omega_\infty)$, and $\mathfrak{a}(([\mathbb{R}]^{<\omega})^\omega_\infty)$, all with respect to the relation $\not\approx $ (strictly speaking, the first w.r.t. (\ref{f|g}) and the latter three w.r.t. \thref{f|g2}), both for their intrinsic interest and their relation to $\mathfrak{a}(\top^*)$.

The first thing to note is that it makes no difference if we replace $\mathbb{R}$ in each of these invariants with $\mathbb{Q}$ because every real can be approximated arbitrarily closely by rationals.  From a forcing point of view, it is nicer to deal with $\mathbb{Q}$ as rationals are countable and absolute w.r.t. transitive models of ZFC, in contrast to reals.

The next thing to note is that all these invariants are at least as large as $\mathfrak{b}$, the minimum cardinality of a $\leq^*$-unbounded subset of $\omega^\omega$ (where $f\leq^*g\Leftrightarrow\forall^\infty n\in\omega(f(n)\leq^*g(n))$).  We also immediately obtain the inequalities $\mathfrak{a}(([\mathbb{R}]^{<\omega}\backslash\{0\})^\omega)\leq\mathfrak{a}(\mathbb{R}^\omega)$ and $\mathfrak{a}(([\mathbb{R}]^{<\omega})^\omega_\infty)\leq\mathfrak{a}(([\mathbb{R}]^{\leq1})^\omega_\infty)$.  We show later that these inequalities can consistently be strict, specifically we show in \thref{random} that they differ in the random model.

It is not difficult to see that $\mathfrak{a}=\mathfrak{a}(([\omega]^{<\omega})^\omega_\infty)$.  It seems possible that also $\mathfrak{a}=\mathfrak{a}(([\mathbb{R}]^{<\omega})^\omega_\infty)$, although this is not investigated in this paper.

\subsection{Main Results and Open Questions}

Our main result is that the cardinal invariants just defined can be consistently less than $\mathfrak{c}$, specifically that this holds in the Sacks model.  This is done for functions in \thref{fSacks} and for projections in \thref{PSacks}.  However, it remains open whether we can actually prove $\mathfrak{a}(\top^*)=\aleph_1$ in ZFC.  Indeed, while we have $\mathfrak{a}(\mathbb{R}^\omega)\geq\mathfrak{b}$, for example, and hence consistently $\mathfrak{a}(\mathbb{R}^\omega)(\geq\mathfrak{b})>\aleph_1$, it is not clear that the same is true for bounded functions, i.e. whether we consistently have $\mathfrak{a}([0,1]^\omega)>\aleph_1$.  There is even evidence to the contrary.  Specifically, we can show that a weaker (i.e. smaller) version of $\mathfrak{a}([0,1]^\omega)$ is indeed equal to $\aleph_1$ in ZFC.
\begin{dfn}
Let $\kappa$ be a cardinal, $S$ a set and $R$ a (symmetric) binary relation on $S$.  We define \[\mathfrak{p}_\kappa((S,R))=\inf\{|A|:A\subseteq S\wedge\forall s\in S\exists a\in A\neg(sRa)\wedge\forall B\in[A]^{<\kappa}\exists s\in S\forall b\in B(sRb)\}.\]
\end{dfn}
Note that we always have $\mathfrak{p}_\kappa((S,R))\leq\mathfrak{a}_\kappa((S,R))$.  Also $\mathfrak{p}=\mathfrak{p}_{\aleph_0}(([\omega]^\omega,\not|^*))$, where $\mathfrak{p}$ is the standard pseudointersection number.

\begin{thm}
In ZFC we have $\mathfrak{p}([0,1]^\omega)=\aleph_1$.
\end{thm}

\paragraph{Proof:} First we claim that we have a net $(X_x)_{x\in[\aleph_1]^{<\omega}}\subseteq[\omega]^\omega$ satisfying the following properties
\begin{eqnarray}
X_x &\subseteq& X_{x\backslash\{\max(x)\}},\textrm{ for all }x\in[\aleph_1]^{<\omega},\label{X1}\\
\textrm{and }X_x\cap X_y &=& 0,\textrm{ for all }x,y\in[\aleph_1]^{<\omega}\textrm{ with }\max(x)=\max(y)\textrm{ and }x\neq y.\label{X2}
\end{eqnarray}
To see this, recursively define $(X_x)_{x\in[\aleph_1]^{<\omega}}$ as follows.  Set $X_0=\omega$ and, once $(X_x)_{x\in[\xi]^{<\omega}}$ has been defined for some $\xi<\aleph_1$, note that $[\xi]^{<\omega}$ is countable and hence we have a sequence $(x_n)\subseteq[\xi]^{<\omega}$ such that, for each $x\in[\xi]^{<\omega}$, there exists infinitely many $n\in\omega$ such that $x=x_n$.  For each $n\in\omega$, recursively pick $m_n\in X_{x_n}\backslash\{m_k:k<n\}$.  Once this is done set $X_{x\cup\{\xi\}}=\{m_n:x=x_n\}$, for each $x\in[\xi]^{<\omega}$.


Define $(f_\xi)_{\xi<\aleph_1}\subseteq\mathbb{R}^{\omega\times\omega}$ as follows.  For each $x\in[\aleph_1]^{<\omega}$, $n\in X_x$ and $m\in\omega$ let us set $f_{\max(x)}(n,m)=|x|/(m+1)$, noting that (\ref{X2}) ensures that each $f_\xi$ is a well defined function on a subset of $\omega$ (we may define $f_\xi$ arbitrarily for the other values).

Now take $f\in[0,1]^{\omega\times\omega}$ and assume that $f\not\approx f_\xi$ for all $\xi\in\aleph_1$.  By the pigeonhole principle, we must therefore have $k\in\omega$, $F\in[\omega\times\omega]^{<\omega}$ and $X\in[\aleph_1]^{\aleph_1}$ such that $|f(n,m)-f_\xi(n,m)|>1/k$, for all $\xi\in X$ and $(n,m)\in\omega\times\omega\backslash F$.  But then, letting $x$ be any $k$-element subset of $X$ and letting $n\in X_x$ be large enough that $(n,k)\notin F$, we see that $f_{\xi_j}(n,k)=(j+1)/(k+1)$, for all $j\in k$, where $(\xi_j)_{j\in k}$ is the increasing enumeration of $x$.  As $f(n,k)\in[0,1]$, we must therefore have $|f(n,k)-f_\xi(n,k)|\leq1/(k+1)$ for some $\xi\in x\subseteq X$, a contradiction. $\Box$\\

If we go in between bounded and unbounded functions, i.e. if we look at $\prod I_n$ for intervals $(I_n)\subseteq\mathbb{R}$ which increase in size sufficiently fast, then we see that consistently $\mathfrak{p}(\prod I_n)>\aleph_1$ by the same argument used to show $\mathfrak{a}(\mathbb{R}^\omega)=\mathfrak{c}$ in \thref{random}.  On the other hand, it is open whether a slight variant of $\mathfrak{p}([0,1]^\omega)$, namely $\mathfrak{p}(([0,1]^{\leq1})^\omega_\infty)$, is still $\aleph_1$ in ZFC.  It is not clear which versions of $\mathfrak{p}$-like and $\mathfrak{a}$-like invariants for functions are closest to those for projections, but an answer to this last question for $\mathfrak{p}(([0,1]^{\leq1})^\omega_\infty)$ will most likely shed light on whether $\mathfrak{p}(\top^*)=\aleph_1$ in ZFC.  If $\mathfrak{p}(\top^*)$ is consisently $>\aleph_1$ then the same is true for $\mathfrak{a}(\top^*)$, while a proof that $\mathfrak{p}(\top^*)=\aleph_1$ in ZFC may well provide a way of also proving the stronger statement that $\mathfrak{a}(\top^*)=\aleph_1$ in ZFC.

\section{Consistency Results}

\subsection{Functions on $\omega$}

\begin{dfn}[\cite{d} Definition 6.3.1]
A forcing notion $\mathbb{P}$ is $\omega^\omega$-bounding if, for all $p\in\mathbb{P}$ and names $\dot{f}$ for elements of $\omega^\omega$, there exists $q\leq p$ and $g\in\omega^\omega$ such that $q\Vdash\forall n\in\omega(\dot{f}(n)\leq g(n))$.
\end{dfn}

A Cohen indestructible infinite MAD family of $[\omega]^\omega$ is constructed in \cite{e} VIII Theorem 2.3.  Our approach, using instead the $\omega^\omega$-bounding property (which does not apply to Cohen forcing), is more similar to \cite{f} Lemma III.1.

\begin{thm}\thlabel{bounding}
If CH holds in the ground model $V$ and $\mathbb{P}$ is an $\omega^\omega$-bounding proper forcing notion with $(|\mathbb{P}|=\aleph_1)^V$, we have lim-inf disjoint $\mathcal{F}\subseteq([\mathbb{Q}]^{<\omega}\backslash\{0\})^\omega$ in $V$ such that $\mathds{1}\Vdash_\mathbb{P}\mathcal{F}$ is maximal.
\end{thm}

\paragraph{Proof:} Thanks to CH and properness, there are only $\aleph_1$ many nice names for elements of $([\mathbb{Q}]^{\leq1})^\omega_\infty$.  Let $(p_\xi,\tau_\xi)_{\xi\in\omega_1}$ enumerate all pairs of elements of $\mathbb{P}$ and such nice names.  We construct $(f_\xi)_{\xi\in\omega_1}\subseteq([\mathbb{Q}]^{<\omega}\backslash\{0\})^\omega$ by recursion as follows.  Say we have constructed $(f_\xi)_{\xi\in\gamma}$.  If
\begin{equation*}
p_\gamma\not\Vdash\forall\xi\in\gamma(f_\xi\not\approx \tau_\gamma)
\end{equation*}
then set $f_\gamma$ to be, say, $f_0$ (or simply leave $f_\gamma$ undefined).  Otherwise, let $(\xi_n)_{n\in\omega}$ be an enumeration of all $\xi\in\gamma$.  Thanks to the $\omega^\omega$-bounding property, there exists $(N_n),(m_n)\subseteq\omega$ in the ground model and $p'_\gamma\leq p_\gamma$ such that
\begin{equation}\label{b1}
p'_\gamma\Vdash\forall n\in\omega\forall k\geq m_n\forall s\in f_{\xi_n}(k)\forall t\in\tau_\gamma(k)(|s-t|>1/N_n).
\end{equation}
Again, by the $\omega^\omega$-bounding property, we may find $p''_\gamma\leq p'_\gamma$ and $f\in([\mathbb{Q}]^{<\omega})^\omega$ such that we have $p''_\gamma\Vdash\forall n\in\omega(\tau_\gamma(n)\subseteq f(n))$.  For all $n\in\omega$ let \[f_\gamma(n)=\{r\in\mathbb{Q}:\exists p\leq p''_\gamma(p\Vdash r\in\tau_\gamma(n))\}\subseteq f(n).\]  Note that (\ref{b1}) implies that $f_\xi\not\approx f_\gamma$, for all $\xi<\gamma$.  As $\gamma$ is a countable ordinal, we may also recursively make $f_\gamma(n)\neq0$, for all $n\in\omega$, while still having $f_\xi\not\approx f_\gamma$, for all $\xi<\gamma$.  This completes the recursion.

Let $\mathcal{F}=\{f_\xi:\xi\in\omega_1\}$. If we had $\mathds{1}\not\Vdash_\mathbb{P}\mathcal{F}$ is maximal, then there would exist $p\in\mathbb{P}$ and a nice name $\tau$ for an element of $([\mathbb{Q}]^{\leq1})^\omega_\infty$ such that $p\Vdash\forall\xi\in\omega_1(f_\xi\not\approx \tau)$. Taking $\gamma\in\omega_1$ such that $p=p_\gamma$ and $\tau=\tau_\gamma$ we have
\begin{eqnarray*}
p_\gamma &\Vdash& \forall\xi\in\gamma(f_\xi\not\approx \tau)\qquad\textrm{and hence}\\
p''_\gamma &\Vdash& \forall n\in\omega(\tau_\gamma(n)\subseteq f_\gamma(n)),
\end{eqnarray*}
by our recursive construction.  But we also have $p_\gamma\Vdash(f_\gamma\not\approx \tau_\gamma)$ and thus $p''_\gamma\leq p_\gamma$ forces contradictory statements, a contradiction. $\Box$\\

\begin{dfn}[\cite{d} Definition 7.4.8]
Two reals $f,g\in\omega^\omega$ are \emph{eventually different} if $\{n\in\omega:f(n)=g(n)\}$ is finite.
\end{dfn}

\begin{prp}\thlabel{ediff}
If $V$ is a transitive model of ZFC and $f\in\omega^\omega$ is eventually different from every element of $\omega^\omega\cap V$ then $f$ is lim-inf disjoint from every element of $\mathbb{R}^\omega\cap V$.
\end{prp}

\paragraph{Proof:} Simply note that, for any $g\in\mathbb{R}^\omega\cap V$, there exists $h\in\omega^\omega\cap V$ such that $|g(n)-h(n)|\leq1/2$, for all $n\in\omega$, and hence $\{n\in\omega:|f(n)-g(n)|<1/2\}\subseteq\{n\in\omega:|f(n)-h(n)|<1\}$ is finite. $\Box$\\

\begin{dfn}
For any index set $I$, let $\mu$ be the standard product measure on $2^I$ and let $\mathbb{B}_I\in V$ be the collection of (equivalence classes of) Baire subsets of $2^I$ ordered by inclusion modulo null subsets.
\end{dfn}

\begin{cor}\thlabel{random}
If $V$ is a c.t.m. of ZFC+CH and $(\kappa=\kappa^\omega)^V$ then in $V^{\mathbb{B}_I}$ we have
\begin{eqnarray}
\aleph_1 &=& \mathfrak{a}(([\mathbb{R}]^{<\omega}\backslash\{0\})^\omega)=\mathfrak{a}(([\mathbb{R}]^{<\omega})^\omega_\infty)\quad\textrm{and}\label{r1}\\
\mathfrak{c} &=& \mathfrak{a}(\mathbb{R}^\omega)=\mathfrak{a}(([\mathbb{R}]^{\leq1})^\omega_\infty)=\kappa.\label{r2}
\end{eqnarray}
\end{cor}

\paragraph{Proof:} Any new real in $V^{\mathbb{B}_I}$ will be in $V^{\mathbb{B}_J}$ for some countable $J\subseteq I$.  But $V$ and $\mathbb{B}_J$ satisfy the hypotheses of \thref{bounding} so (\ref{r1}) follows from this.  On the other hand, any collection of $<\kappa$ reals will be in a proper submodel of $V^{\mathbb{B}_I}$, over which $V^{\mathbb{B}_I}$ will contain an eventually different real, so (\ref{r2}) follows from \thref{ediff}. $\Box$\\

\begin{dfn}[\cite{d} Definition 6.3.38]
A forcing notion $\mathbb{P}$ has the Sacks property if, for all $p\in\mathbb{P}$ and names $\dot{f}$ for elements of $\omega^\omega$, there exists $q\leq p$ and $F\in\mathscr{P}(\omega)^\omega$ such that $|F(n)|=2^n$, for all $n\in\omega$, and $q\Vdash\forall n\in\omega(\dot{f}(n)\in F(n))$.
\end{dfn}

\begin{thm}\thlabel{sacks}
If CH holds in the ground model $V$ and $\mathbb{P}$ is a proper forcing notion with the Sacks property and $(|\mathbb{P}|=\aleph_1)^V$ then we have lim-inf disjoint $\mathcal{F}\subseteq\mathbb{Q}^\omega$ in $V$ such that $\mathds{1}\Vdash_\mathbb{P}\mathcal{F}$ is maximal.
\end{thm}

\paragraph{Proof:} It is immediate that any forcing notion with the Sacks property is $\omega^\omega$-bounding so we may proceed as in the proof of \thref{bounding} up to and including (\ref{b1}).  Yet again, by $\omega^\omega$-bounding, we may decrease $p'_\gamma$ and increase $(m_n)$ if necessary so that \[p'_\gamma\Vdash\forall n\in\omega(|\{k\in m_{n+1}\backslash m_n:\tau_\gamma(k)\neq0\}|\geq2^n).\]  Then, by the Sacks property, there exists $p''_\gamma\leq p'_\gamma$ and a function $F$ on $\omega$ such that, for each $n\in\omega$, $F(n)$ is a non-empty collection of functions from $m_{n+1}\backslash m_n$ to $[\mathbb{Q}]^{\leq1}$ of size at most $2^n$ satisfying
\begin{eqnarray}
&& p''_\gamma\Vdash\forall n(\tau_\gamma\upharpoonright m_{n+1}\backslash m_n\in F(n))\qquad\textrm{and}\label{cond2}\\
&& \forall n\forall f\in F(n)\exists p\leq p''_\gamma(p\Vdash\tau_\gamma\upharpoonright m_{n+1}\backslash m_n=f).\label{cond3}
\end{eqnarray}
Let $f_\gamma\in([\mathbb{Q}]^{\leq1})^\omega$ satisfy
\begin{eqnarray}
&& \forall n\forall f\in F(n)\exists k\in m_{n+1}\backslash m_n({f_\gamma(k)}={f(k)}\neq0)\qquad\textrm{and}\label{cond4}\\
&& \forall n\forall k\in m_{n+1}\backslash m_n\exists f\in F(n)(f_\gamma(k)=f(k)).\label{cond5}
\end{eqnarray}
(\ref{b1}), (\ref{cond3}) and (\ref{cond5}) imply that $f_\xi\not\approx f_\gamma$, for all $\xi<\gamma$.  As $\gamma$ is a countable ordinal, we may again recursively make $|f_\gamma(n)|=1$, for all $n\in\omega$, while still having $f_\xi\not\approx f_\gamma$, for all $\xi<\gamma$.  On the other hand, (\ref{cond2}) and (\ref{cond4}) imply that
\begin{eqnarray*}
p''_\gamma &\Vdash& \forall n\exists k\in m_{n+1}\backslash m_n(\tau_\gamma(k)=f_\gamma(k)\neq0)\\
&\Vdash& \exists^\infty n(\tau_\gamma(n)=f_\gamma(n)\neq0),
\end{eqnarray*}
leading to essentially the same contradiction as in the proof of \thref{bounding}. $\Box$\\

Let $\mathbb{S}$ be Sacks forcing and, for any ordinal $\xi$, let $\mathbb{S}_\xi$ its $\xi$-step countable support iteration.

\begin{cor}\thlabel{fSacks}
If CH holds in the ground model $V$ then in $V^{\mathbb{S}_{\omega_2}}$ we have \[\mathfrak{a}(\mathbb{R}^\omega)=\mathfrak{a}(([\mathbb{R}]^{\leq1})^\omega_\infty)=\aleph_1(<\aleph_2=\mathfrak{c}).\]
\end{cor}

\paragraph{Proof:} As shown in \cite{f} Theorem III.2, it suffices to find $\mathcal{F}\in V\cap\mathbb{Q}^\omega$ such that $\mathds{1}\Vdash_{\mathbb{S}_{\omega_1}}\mathcal{F}$ is maximal.  As $\mathbb{S}$ is proper and has the Sacks property, the same applies to its iterations and hence, as $\mathbb{S}_{\omega_1}$ has a dense subset of size $\aleph_1$, this follows from \thref{sacks}. $\Box$\\

\subsection{Projections}

As before, let $H$ be $l^2$, a separable infinite dimensional Hilbert space with (Hilbert) basis $(e_n)$.

\begin{lem}\thlabel{prolem}
Say we have orthogonal unit vectors $v_0,\ldots,v_{n-1}$ in $H$ and a projection $P$ on $H$ such that $v_i\perp Pv_j$, for all distinct $i,j\in n$.  Then $||PP_{\mathrm{span}_{i\in n}(v_i)}||=\max_{i\in n}||Pv_i||$.
\end{lem}

\paragraph{Proof:} For all distinct $i,j\in n$, we have $v_i\perp Pv_j$ and hence $\langle Pv_i, Pv_j\rangle=\langle v_i, Pv_j\rangle=0$, i.e. $Pv_i\perp Pv_j$.  Thus, for all $\alpha_0,\ldots,\alpha_{n-1}\in\mathbb{F}$ such that $|\alpha_0|^2+\ldots+|\alpha_{n-1}|^2=1$,
\begin{eqnarray*}
||P(\alpha_0v_0+\ldots+\alpha_{n-1}v_{n-1})||^2 &=& ||\alpha_0Pv_0+\ldots+\alpha_{n-1}Pv_{n-1}||^2\\
&=& |\alpha_0|^2||Pv_0||^2+\ldots+|\alpha_{n-1}|^2||Pv_{n-1}||^2\\
&\leq& \max_{i\in n}||Pv_i||^2.
\end{eqnarray*}
i.e. $||PP_{\mathrm{span}_{i\in n}(v_i)}||\leq\max_{i\in n}||Pv_i||$.  The reverse inequality is immediate. $\Box$\\

\begin{dfn}
Let $(I_n)$ be a sequence of finite intervals of $\omega$, i.e. of sets of the form $m\backslash n$ for some $n\in m\in\omega$.  We say $(I_n)$ is an \emph{interval subpartition} (of $\omega$) if $\max(I_n)<\min(I_{n+1})$, for all $n\in\omega$.  $(I_n)$ is an \emph{interval partition} if $\min(I_0)=0$ and $\max(I_n)+1=\min(I_{n+1})$, for all $n\in\omega$.
\end{dfn}

\begin{dfn}
$V\subseteq H$ is a \emph{block subspace} (w.r.t. $(e_n)$) if there exists an interval (sub)partition $(I_n)$ and (unit) vectors $(v_n)\subseteq H$ such that $v_n\in\mathrm{span}_{i\in I_n}(e_i)$, for all $n\in\omega$, and $V=\overline{\mathrm{span}}(v_n)$.  $V\subseteq H$ is a \emph{generalized block subspace} (w.r.t. $(e_n)$) if there exists an interval (sub)partition $(I_n)$ and finite subspaces $(F_n)\subseteq H$ such that $F_n\subseteq\mathrm{span}_{i\in I_n}(e_i)$, for all $n\in\omega$, and $V=\bigoplus F_n(=\overline{\sum F_n}$, i.e. the closure of the subspace of finite linear combinations of elements of $\bigcup F_n)$.
\end{dfn}

Let $\mathbb{G}$ be an absolute (w.r.t. countable transitive models of set theory) countable dense subfield of the scalar field $\mathbb{F}$ that is closed under taking square roots (eg. the algebraic numbers).  Then we define $\mathbb{G}$-block subspaces analogously, with the extra requirement that each $v_n$ is a $\mathbb{G}$-vector, i.e. in the $\mathbb{G}$-span of the basis vectors $(e_i)$.

Note that the projections onto infinite dimensional $\mathbb{G}$-block subspaces are $\leq^*$-dense among all infinite rank projections.  For any $P,Q,R\in\mathcal{P}(H)$ we have \[R\leq^*Q\Rightarrow||\pi(PR)||=||\pi(PQR)||\leq||\pi(PQ)||\] and hence, to verify that some almost disjoint family $\mathcal{P}$ of infinite rank projections is maximal we need only verify that, for all projections $R$ onto infinite dimensional $\mathbb{G}$-block subspaces of $H$, there exists $P\in\mathcal{P}$ such that $||\pi(PR)||=1$.

\begin{thm}\thlabel{sacksp}
If CH holds in the ground model $V$ and $\mathbb{P}$ is a proper forcing notion with the Sacks property such that $(|\mathbb{P}|=\aleph_1)^V$ then there exists an almost disjoint family $\mathcal{P}\subseteq\mathcal{P}_\infty(H)$ in $V$ such that $\mathds{1}\Vdash_\mathbb{P}\mathcal{P}$ is maximal.
\end{thm}

\paragraph{Proof:} Thanks to CH and properness, there are only $\aleph_1$ many nice names for projections onto infinite dimensional $\mathbb{G}$-block subspaces of $H$.  In fact, all bounded linear operators can be considered as reals by looking at their matrix representations.  Indeed, if we are being truly formal, we have to deal with these matrices as coding the projections we are talking about, seeing as $l^2$, the domain of the projections, becomes larger in any extension adding reals.  That aside, let $(p_\xi,\tau_\xi)_{\xi\in\omega_1}$ enumerate all pairs of elements of $\mathbb{P}$ and such nice names.  We construct $(P_\xi)_{\xi\in\omega_1}\subseteq\mathcal{P}$ and interval subpartitions $((K^\xi_n))_{\xi\in\omega_1}$ (each $(K^\xi_n)$ corresponding to the generalized blocks of $P_\xi$) by recursion as follows.  Say we have constructed $(P_\xi)_{\xi\in\gamma}$.  If
\begin{equation*}
p_\gamma\not\Vdash\forall\xi\in\gamma||\pi(\tau_\gamma P_\xi)||<1
\end{equation*}
then set $P_\gamma$ to be, say, $P_0$ (or simply leave $P_\gamma$ undefined).  Otherwise, let $(\xi_n)_{n\in\omega}$ be an enumeration of all $\xi\in\gamma$ and let $(\dot{I}_n)$ and $(\dot{v}_n)$ be names for the interval partition and corresponding unit vectors that define the block subspace $\mathcal{R}(\tau_\gamma)$.  Thanks to the Sacks property which in turn implies the $\omega^\omega$-bounding property, there exists $(N_n)\subseteq\omega$ in the ground model and $p'_\gamma\leq p_\gamma$ such that
\begin{equation*}
p'_\gamma\Vdash\forall n||\pi(\tau_\gamma P_{\xi_n})||<N_n/(N_n+1).
\end{equation*}
Then, again thanks to the $\omega^\omega$-bounding property, there exists increasing $(m_n)\subseteq\omega$, with $m_0=0$, in the ground model $V$ and $p''_\gamma\leq p'_\gamma$ such that
\begin{equation}\label{cond1p}
p''_\gamma\Vdash\forall l<n||P_{\xi_l}P_{\mathrm{span}\{\dot
{v}_k:k\geq m_n\}}||<N_l/(N_l+1).
\end{equation}
Yet again thanks to the $\omega^\omega$-bounding property, there exists increasing an interval partition $(J_n)$ of $\omega$ in the ground model $V$ and $p'''_\gamma\leq p''_\gamma$ such that
\begin{equation*}
p'''_\gamma\Vdash\forall n\exists m(\dot{I}_m\subseteq J_n).
\end{equation*}
Let $(K^\gamma_n)$ be an interval subpartition such that, for all $n\in\omega$, \[\max\{k:K^\gamma_n\cap K^{\xi_l}_k\neq0\}<\min\{k:K^\gamma_{n+1}\cap K^{\xi_l}_k\neq0\},\] for all $l<n$, $K^\gamma_n$ contains only intervals $J_k$ such that $k\geq m_n$ and contains more than $(n+1)(2^n-1)$ such intervals, and hence
\begin{equation*}
p'''_\gamma\Vdash\forall n((\dot{I}_k\subseteq K^\gamma_n\Rightarrow k\geq m_n)\wedge|\{k:\dot{I}_k\subseteq K^\gamma_n\}|>(n+1)(2^n-1)).
\end{equation*}
By the Sacks property, there exists $p''''_\gamma\leq p'''_\gamma$ and function $F$ on $\omega$ such that, for each $n\in\omega$, $|F(n)|\leq2^n$ and each element of $F(n)$ is a collection of $\mathbb{G}$-vectors such that
\begin{eqnarray}
&& p''''_\gamma\Vdash\forall n(\{\dot{v}_k:\dot{I}_k\subseteq K^\gamma_n\}\in F(n))\label{cond2p}\\
&\textrm{and}& \forall n\forall \mathcal{V}\in F(n)\exists p\leq p''_\gamma(p\Vdash\{\dot{v}_k:\dot{I}_k\subseteq K^\gamma_n\}=\mathcal{V}).\label{cond3p}
\end{eqnarray}

For each $n$, let $\mathcal{V}_0,\ldots,\mathcal{V}_{|F(n)|-1}$ enumerate the elements of $F(n)$ and define vectors $u_n^0,\ldots,u_n^{|F(n)|-1}$ recursively as follows.  Let $u_n^0$ be any element of $\mathcal{V}_0$.  Once $u_n^j$ has been defined, for $j<i$, let $u_n^i$ be any unit ($\mathbb{G}$-)vector in $\mathrm{span}(\mathcal{V}_i)$ that is perpendicular to $u_0,\ldots,u_{i-1}$ and $P_ku_j$, for all $k\in n$ and $j\in i$.  As $\dim\mathrm{span}(\mathcal{V})=|\mathcal{V}|>(n+1)(2^n-1)$, for all $\mathcal{V}\in F(n)$, and $|F(n)|\leq2^n$ we may continue this recursion until we have defined $u_n^{|F(n)|-1}$.  Let $U_n=\mathrm{span}(u_n^i)_{i\in|F(n)|}$ and let $P_\gamma$ be the projection onto $\bigoplus U_n$.  \thref{prolem}, (\ref{cond1p}) and (\ref{cond3p}) imply that \[\forall l\in\omega(||P_{\xi_l}P_{\bigoplus_{n>l}U_n}||<N_l/(N_l+1))\] and hence that $\forall n(||\pi(P_{\xi_n}P_\gamma)||<N_n/(N_n+1))$, as $\pi(P_{\bigoplus_{n>l}U_n})=\pi(P_\gamma)$, which, in particular, means that $P_\gamma$ is almost disjoint from $(P_\xi)_{\xi\in\gamma}$.  On the other hand, (\ref{cond2p}) implies that
\begin{equation*}
p''''_\gamma\Vdash\forall n(\mathcal{R}(\tau_\gamma)\cap\mathcal{R}(P_\gamma)\cap\mathrm{span}_{i\in K_n}(e_i)\neq0)
\end{equation*}
which, in particular, means that $p''''_\gamma\Vdash\mathcal{R}(\tau_\gamma)\cap\mathcal{R}(P_\gamma)$ is infinite dimensional, so
\begin{equation}
p''''_\gamma\Vdash||\pi(\tau_\gamma P_\gamma)||=1.
\end{equation}
This completes the recursive construction.

If we had $\mathds{1}\not\Vdash_\mathbb{P}\mathcal{P}$ is maximal, then there would exist $p\in\mathbb{P}$ and a nice name $\tau$ for a projection onto an infinite dimensional $\mathbb{G}$-block subspace of $H$ such that \[p\Vdash\forall\xi\in\omega_1(||\pi(\tau_\xi P_\xi)||<1).\] Taking $\gamma\in\omega_1$ such that $p=p_\gamma$ and $\tau=\tau_\gamma$ we have
\begin{eqnarray*}
p_\gamma &\Vdash& \forall\xi\in\gamma(||\pi(\tau_\gamma P_\xi)||<1)\qquad\textrm{and hence}\\
p''''_\gamma &\Vdash& (||\pi(\tau_\gamma P_\gamma)||=1),
\end{eqnarray*}
by our recursive construction.  But we also have \[p_\gamma\Vdash||\pi(\tau_\gamma P_\gamma)||<1\] and thus $p''''_\gamma\leq p_\gamma$ forces contradictory statements, a contradiction. $\Box$\\

\begin{cor}
If CH holds in the ground model $V$ then there exists almost disjoint $\mathcal{P}\subseteq\mathcal{P}_\infty(H)$ in $V$ such that $\mathds{1}\Vdash_{\mathbb{S}_{\omega_2}}\mathcal{P}$ is maximal.
\end{cor}

\paragraph{Proof:} As shown in \cite{f} Theorem III.2, it suffices to find $\mathcal{P}$ such that $\mathds{1}\Vdash_{\mathbb{S}_{\omega_1}}\mathcal{P}$ is maximal.  As $\mathbb{S}$ is proper and has the Sacks property, the same applies to its iterations and hence, as $\mathbb{S}_{\omega_1}$ has a dense subset of size $\aleph_1$, this follows from \thref{sacksp}. $\Box$\\

\begin{cor}\thlabel{PSacks}
If CH holds in the ground model $V$ then $V^{\mathbb{S}_{\omega_2}}$ satsifies $\mathfrak{a}(\top^*)=\aleph_1(<\aleph_2=\mathfrak{c})$.
\end{cor}

\end{document}